\documentclass[a4paper,leqno,11pt]{amsart}
\usepackage{amsfonts,amssymb,verbatim,amsmath,amsthm,latexsym,textcomp,amscd}
\usepackage{latexsym,amsfonts,amssymb,epsfig,verbatim}
\usepackage{amsmath,amsthm,amssymb,latexsym,graphics,textcomp}
\usepackage{graphicx}
\usepackage{color}
\usepackage{url}
\usepackage{enumerate}
\usepackage[mathscr]{euscript}
\usepackage{tikz}
\usepackage{dsfont}
\usetikzlibrary{matrix}
\usepackage{sseq}
\usepackage{hyperref}

\input xy
\xyoption{all}

\theoremstyle{plain}
\newtheorem*{thma}{Theorem A}
\newtheorem*{thmb}{Theorem B}
\newtheorem*{thmc}{Theorem C}

\theoremstyle{definition}
\newtheorem{theorem}{Theorem}[section]
\newtheorem{prop}[theorem]{Proposition}
\newtheorem{lemma}[theorem]{Lemma}

\newtheorem{rmk}[theorem]{Remark}
\newtheorem{subsec}[theorem]{}

\newtheorem{thm}[theorem]{Theorem}
\newtheorem{notation}[theorem]{Notation}

\theoremstyle{remark}
\swapnumbers

\newcommand{\C}{{\mathbb C}}

\newcommand{\Z}{{\mathbb{Z}}}
\newcommand{\R}{{\mathbb R}}
\newcommand{\Q}{{\mathbb Q}}

\newcommand\DD{{\mathcal D}}

\newcommand\FF{{\mathcal F}}

\newcommand\LL{{\mathcal L}}
\newcommand\MM{{\mathcal M}}

\newcommand\PP{{\mathcal P}}

\newcommand{\Sa}{\mathcal S}

\newcommand\PMF{{\PP\kern-2pt\MM\FF}}

\newcommand\PML{{\PP\kern-2pt\MM\LL}}

\newcommand{\fsubd}{\mathrel{{\scriptstyle\searrow}\kern-1ex^d\kern0.5ex}}
\newcommand{\bsubd}{\mathrel{{\scriptstyle\swarrow}\kern-1.6ex^d\kern0.8ex}}
\newcommand{\fsubeq}{\mathrel{\raise-.7ex\hbox{$\overset{\searrow}{=}$}}}
\newcommand{\bsubeq}{\mathrel{\raise-.7ex\hbox{$\overset{\swarrow}{=}$}}}

\newcommand{\tsh}[1]{\left\{\kern-.9ex\left\{#1\right\}\kern-.9ex\right\}}

\newenvironment{myeq}[1][]
{\stepcounter{theorem}\begin{equation}\tag{\thetheorem}{#1}}
{\end{equation}}

\newenvironment{mysubsection}[2][]
{\begin{subsec}\begin{upshape}\begin{bfseries}{#2.}
			\end{bfseries}{#1}}
		{\end{upshape}\end{subsec}}

\makeatletter
\@namedef{subjclassname@2020}{%
  \textup{2020} Mathematics Subject Classification}
\makeatother

\title{$p$-local decompositions of projective Stiefel manifolds}
\author{Samik Basu, Debanil Dasgupta, Shilpa Gondhali, Swagata Sarkar}

\address{Stat-Math Unit,
Indian Statistical Institute,
B. T. Road, Kolkata-700108, India.}
\email{samik.basu2@gmail.com; samikbasu@isical.ac.in}

\address{Stat-Math Unit,
Indian Statistical Institute,
B. T. Road, Kolkata-700108, India.}
\email{debanil12@gmail.com;}

\address{Department of Mathematics,
BITS Pilani K K Birla Goa Campus,
NH 17B, Bypass Road, Zuarinagar, 
Sancoale, Goa 403726. India.}
\email{shilpag@goa.bits-pilani.ac.in; shilpa.s.gondhali@gmail.com}

\address{School of Mathematical Sciences,  
UM-DAE Centre for Excellence in Basic Sciences, 
University of Mumbai, Vidyanagari Campus, Kalina, 
Santacruz (East),  Mumbai - 400098,   India.}
\email{ swagata.sarkar@cbs.ac.in}


\begin{document}

\maketitle

\begin{abstract}
The main objective of this paper is to analyze the $p$-local homotopy type of the complex projective Stiefel manifolds, and other analogous quotients of Stiefel manifolds. We take the cue from a result of Yamaguchi about the $p$-regularity of the complex Stiefel manifolds which lays down some hypotheses under which the Stiefel manifold is $p$-locally a product of odd dimensional spheres. We show that in many cases, the projective Stiefel manifolds are $p$-locally a product of a complex projective space and some odd dimensional spheres. As an application, we prove that in these cases, the $p$-regularity result of Yamaguchi is also $S^1$-equivariant. 
\end{abstract}

\section{Introduction}
The complex Stiefel manifold $W_{n,k}$ is the space of $k$-orthonormal vectors in $\C^n$. This has an $S^1$-action by multiplying each vector by a complex number of norm $1$. This paper attempts to prove homotopical results about the quotients of the Stiefel manifolds under this action and its variants. The cohomology of the projective Stiefel manifold, defined as $PW_{n,k}:=W_{n,k}/S^1$, was computed in \cite{AGMP}, and this was used to reveal conditions for equivariant maps \cite{PP13}. 

The cohomology of the Stiefel manifold is an exterior algebra, so a natural question is how much homotopically alike it is to a product of spheres. If it is so, the fibration $W_{n,k}\to S^{2n-1}$ must have a section $S^{2n-1} \to W_{n,k}$, the existence of which has a precise answer \cite{AW65}, \cite{AT60}. On the other hand, one may try to write down conditions under which $W_{n,k}$ is $p$-regular, that is, when does its $p$-localization become equivalent to a product of spheres? This has been studied by Yamaguchi \cite{Yam89}. His results imply that if $p>n$, $W_{n,k}$ is $p$-regular.  We consider the analogous question for the quotient $PW_{n,k}$.  

The cohomology of $PW_{n,k}$  matches that of $\C P^{n-k}\times \prod_{j=n-k+2}^{n} S^{2j-1}$ if $p>n$. We explore whether, after localization at $p$, the two spaces are homotopy equivalent. Elementary arguments imply that if $p$ is more than half the dimension of $PW_{n,k}$, this result holds. More precisely, (see Theorem \ref{largepdec})
\begin{thma}
If $p>  \frac{2nk-k^2-1}{2}+k-n$,  
\[ \Big(P W_{n,k}\Big)_{(p)} \simeq \Big[ \C P^{n-k} \times S^{2n-2k+3} \times \cdots \times S^{2n-1} \Big]_{(p)}.\]
\end{thma} 

For a significantly better bound on $p$ in the equivalence above, we look at the obstructions to forming a map from $PW_{n,k}$ to the product $\C P^{n-k}\times \prod_{j=n-k+2}^{n} S^{2j-1}$. We see that if $p>n$, these obstructions (localized at $p$) belong to the ``stable range". This leads us to consider the stable homotopy type of the projective Stiefel manifold, for which we prove the following result using certain calculations with the Chern character. (see Theorem \ref{projstsplit}) 
\begin{thmb}
If $p>n$, the projective Stiefel manifold $PW_{n,k}$ stably splits into a wedge of spheres in the $p$-local category. 
\end{thmb}  
The path to proving Theorem B goes via a homotopy theoretic result which says that for a CW-complex of dimension $\leq$ $2p^2 - 2p$, if $p$-local cohomology is torsion-free, and the Chern character takes values in $\Z_{(p)}$, the space stably splits into a wedge of spheres in the $p$-local category (see Theorem \ref{resstsplit}). The stable splitting of the projective Stiefel manifold allows us to construct stable maps to the desired space. The obstruction theory to lift this to the unstable category involves another intricate argument with the Chern character that introduces a new bound on $k$ for which the $p$-local decomposition result holds. (see Theorem \ref{unsplit})
\begin{thmc}
Suppose $p>n+1$ and $k\leq  p+n-\sqrt{p^2+n^2 - 4p +2}$. Then, 
\[ \Big(P W_{n,k}\Big)_{(p)} \simeq \Big[ \C P^{n-k} \times S^{2n-2k+3} \times \cdots \times S^{2n-1} \Big]_{(p)}.\]
\end{thmc}
The homotopical decomposition results for $PW_{n,k}$ imply that with the bound above, the $p$-regularity result for the Stiefel manifold is indeed $S^1$-equivariant (see Theorem \ref{eqstief}). The techniques also imply decomposition results (see Theorem \ref{splitquot}) for the finite cyclic quotients of Stiefel manifolds $W_{n,k;m}$, and that of $P_\ell W_{n,k}$ for $\ell=(l_1,\cdots, l_k)$, $l_i\in \Z$, which is $W_{n,k}/S^1$ with the action on the $i^{th}$ vector by $z^{l_i}$.   

\begin{mysubsection}{Organization}
In \S \ref{cohprojst}, we recall the cohomology of the projective Stiefel manifold and the other quotients of Stiefel manifolds used in the paper. In \S \ref{decvlarge}, we discuss the decomposition results for large primes using elementary arguments. In \S \ref{decst}, the stable decomposition results are proved using the Chern character computations. Finally in \S \ref{unpld}, we complete the $p$-local decomposition results for $p>n$. 
\end{mysubsection}

\begin{notation}
Throughout the document, $p$ stands for an odd prime. We work in the $p$-local category of spaces. For a homogeneous space $X=G/H$, we use the same notation for it's $p$-localization. We write $X_{(p)}$ for the $p$-localization of $X$, and $X_\Q$ for the rationalization. We use $H^\ast(X)$ to denote the cohomology ring $H^\ast(X; \Z_{(p)})$, that is the cohomology ring over $\Z_{(p)}$. The notation $\gamma$ stands for the canonical line bundle over $\C P^\infty$. The following notations are used for homogeneous spaces arising in the document. 
\begin{itemize}
\item $W_{n,k}$ denotes the complex Stiefel manifold $\cong U(n)/U(n-k)$. 
\item $Gr_k(\C^n)$ denotes the complex grassmannian of $k$-planes in $\C^n$.
\item $PW_{n,k}$ denotes the projective Stiefel manifold $W_{n,k}/S^1$ where $S^1$ acts on $W_{n,k}$ via $z\cdot (v_1,\cdots, v_k) = (zv_1,\cdots, zv_k)$. 
\item $Y_{n,k}$ denotes the product $\C P^{n-k}\times \prod_{i=n-k+2}^n S^{2i-1}$. 
\item $W_{n,k;m}$ denotes $W_{n,k}/C_m$, where $C_m$ is the cyclic group of order $m$ acting on $W_{n,k}$ via the natural inclusion $C_m \subset S^1$, and the above action of $S^1$ on $W_{n,k}$. 
\item $L_{m}(2k+1)$ denotes the lens space $S^{2k+1}/C_m=W_{2k+1,1;m}$. 
\item For a tuple $\ell = (l_1,\cdots,l_k)$, $P_\ell W_{n,k}$ denotes the quotient $W_{n,k}/S^1$ where $S^1$ acts on $W_{n,k}$ via $z\cdot (v_1,\cdots, v_k) = (z^{l_1}v_1,\cdots, z^{l_k}v_k)$. 
\end{itemize}
For a tuple $\ell=(l_1,\cdots, l_k)$ and a collection of integers $I=(i_1,\cdots, i_k)$, we write $|I|= \sum_{j=1}^k i_j$ and $l^I=\prod_{j=1}^k l_j^{i_j}$.
\end{notation}


\section{Cohomology of projective Stiefel manifolds} \label{cohprojst}
The purpose of this section is to record the cohomology of projective Stiefel manifolds and other associated quotients of Stiefel manifolds. The cohomology of projective Stiefel manifolds $PW_{n,k}$ with $\Z/p$ coefficients was computed in \cite{AGMP}. An analogous computation works for the cohomology with $\Z_{(p)}$-coefficients, and one has the formula (which is only additive for $p=2$)
\begin{myeq} \label{cohploc} 
H^*(PW_{n,k};\Z_{(p)})\cong \Lambda_{\Z_{(p)}}(\gamma_{n-k+2},\cdots,\gamma_{n})\otimes_{\Z_{(p)}} \frac{\Z_{(p)}[x]}{(\binom{n}{j}x^j \mid n-k+1\leq j\leq n)},
\end{myeq}
with $|\gamma_i|=2i-1$ and $|x|=2$. The computation is carried out using the Serre spectral sequence for the fibration $W_{n,k} \to PW_{n,k} \to \C P^\infty$ \cite[(2.1)]{AGMP}. One identifies  the cohomology of $W_{n,k}$ as the exterior algebra $\Lambda(z_{n-k+1},\cdots, z_n)$ with $|z_j|=2j-1$, and computes the differentials by the fact that the $z_j$ are transgressive, and the equation $d_{2j}(z_j)=\binom{n}{j} x^j$. If $p>n$, one observes that the binomial coefficients $\binom{n}{j}$ are units in $\Z_{(p)}$, so we have the following reduction of \eqref{cohploc}  
\begin{myeq} \label{cohplarge} 
H^*(PW_{n,k};\Z_{(p)})\cong \Lambda_{\Z_{(p)}}(\gamma_{n-k+2},\cdots,\gamma_{n})\otimes \frac{\Z_{(p)}[x]}{(x^{n-k+1})}.
\end{myeq}
In this case, observe that the cohomology of $PW_{n,k}$ matches with that of the product $Y_{n,k} = \C P^{n-k}\times \prod_{i=n-k+2}^n S^{2i-1}$. We now follow \cite[\S 3]{AGMP} to identify the generators $\gamma_j$ of \eqref{cohplarge}. Let $E$ be a contractible space with free $U(n)$-action. We have the following homotopy commutative diagram in which all the squares are homotopy pullbacks,
\begin{myeq}\label{pullback}
\xymatrix{W_{n,k}\ar[r]^-{\pi}\ar[d] &PW_{n,k}\ar[r]^-f\ar[d] &BU(n-k)\ar[d]\\
E\ar[r] &\C P^{\infty}\ar[r]_-{f_0} &BU(n), } 
 \end{myeq}
where $f_0$ classifies the bundle $n\gamma$. This gives rise to the following diagram 
\begin{myeq}\label{PW-H}
\xymatrix{0\ar[d]\\
H^{2j-1}(PW_{n,k};\Z_{(p)})\ar[d]_-{\delta} &0\ar[d]\\
H^{2j}(\C P^{\infty},PW_{n,k};\Z_{(p)})\ar[d]_-{i^*} & H^{2j}(BU(n),BU(n-k);\Z_{(p)})\ar[l]^-{f^*}\ar[d]\\
H^{2j}(\C P^{\infty};\Z_{(p)}) & H^{2j}(BU(n);\Z_{(p)}).\ar[l]^-{f_0^*}}
\end{myeq}
We identify $H^\ast(BU(n),BU(n-k);\Z_{(p)})$ with the ideal of $H^\ast(BU(n);\Z_{(p)}) $ generated by the universal Chern classes $c_j$ for $n-k<j\leq n$, and write $u_j^H=f^*c_j$. In this notation, $\gamma_j^H \in H^{2j-1}(PW_{n,k};\Z_{(p)})$ is defined by the equation
\begin{myeq} \label{choicegen}
 \delta \gamma_j^H = \rho_j= u_j^H - x^{j-(n-k+1)}\frac{\mu_j}{ \mu_{n-k+1}} u_{n-k+1}^H,
\end{myeq}
where $\mu_j=\binom{n}{j}$. In order to observe how this formula makes sense, one should note that $f_0^\ast(c_j)=\binom{n}{j} x^j$. 

\begin{mysubsection}{Other quotients of Stiefel manifolds} 
One may proceed in an analogous manner to the above to write down the cohomology ring structure for the other quotients of Stiefel manifolds. Let $\ell=(l_1,\cdots, l_k)$  such that the $\gcd$ of the $l_i$ is $1$. 
The cohomology of $P_\ell W_{n,k}$ with $\Z/p$-coefficients was computed in \cite{BS17}. The similar calculation with $\Z_{(p)}$-coefficients yields the formula 
$$ H^*(P_\ell W_{n,k};\Z_{(p)}) \cong \Lambda_{\Z_{(p)}}(\gamma_{n-k+2},\cdots, \gamma_n) \otimes \Z_{(p)}[x]/J, $$ 
where $|\gamma_j|= 2j-1$, $|x|=2$ and $J$ is the ideal of $\Z_{(p)}[x]$ generated by the set $\{ \sum_{|I|=j} (-1)^j l^Ix^j\mid n-k<j\leq n \}$.  This formula is obtained by calculating the differentials in the Serre spectral sequence associated to the fibration $W_{n,k} \longrightarrow P_\ell W_{n,k} \longrightarrow \C P^\infty$. It turns out that the exterior algebra generators $z_j$'s of $H^*(W_{n,k};\Z_{(p)})$ are transgressive and $d_{2j}(z_{j})= \sum_{|I|=j} (-1)^j l^Ix^j$.  Note that for a prime $p$ not dividing $\sum_{|I|=n-k+1} l^I$, we have the following reduction 
\begin{myeq} \label{plwnkcoh}
H^*(P_\ell W_{n,k};\Z_{(p)}) \cong \Lambda_{\Z_{(p)}}(\gamma_{n-k+2},\cdots, \gamma_n) \otimes \Z_{(p)}[x]/(x^{n-k+1}).
\end{myeq}
In this case $P_\ell W_{n,k}$ and $Y_{n,k}$ have isomorphic cohomology rings. The analogue of the pullback \eqref{pullback} is the diagram \cite[(2.1)]{BS17}
\begin{myeq}\label{pullbackpl}
\xymatrix{W_{n,k}\ar[r]^-{\pi}\ar[d] & P_\ell W_{n,k}\ar[r]^-{\phi} \ar[d] & Gr_k(\C^n) \ar[d]\\
E\ar[r] &\C P^{\infty}\ar[r]_-{\phi_0} &BU(k), } 
 \end{myeq}
where $\phi_0$ classifies the bundle $\sum \gamma^{l_j}$. One may now consider a diagram similar to \eqref{PW-H} by working with the pair $(BU(k), Gr_k(\C^n))$ to identify the cohomology generators for $P_\ell W_{n,k}$.

The cohomology of $W_{n,k;m}$ with $\Z/p$ coefficients was computed in \cite{GS13}. For $p\nmid m$, this is equivalent to the cohomology of $W_{n,k}$, so the interesting case is when $p\mid m$. Following the same method, the cohomology with $\Z_{(p)}$ coefficients (for $p\mid m$) may be computed. 
For $p>n$, the formula takes the following form 
\begin{myeq} \label{wnkmcoh}
 H^*(W_{n,k;m};\Z_{(p)}) \cong (\Lambda_{\Z_{(p)}}(\gamma_{n-k+1},\gamma_{n-k+2},\cdots, \gamma_n) \otimes \Z_{(p)}[x])/(mx, x^{n-k+1}, \gamma_{n-k+1}x) , 
\end{myeq}
 where $|\gamma_j|= 2j-1$, and  $|x|=2$. The method requires determining the Serre spectral sequence associated to the fibration $S^1 \longrightarrow W_{n,k;m}\longrightarrow PW_{n,k} $ and the only differential $d_2$ sends the degree $1$ class $e$ generating $H^*(S^1;\Z_{(p)})$ to $mx$. Note that the class $e\otimes x^{n-k}$ survives in the $E_\infty$-page detecting the degree $2n-2k+1$ class $\gamma_{n-k+1}$. 
\end{mysubsection}

\section{Decomposition results at very large primes}\label{decvlarge}

In \S \ref{cohprojst}, the expressions for the cohomology of the various quotients of $W_{n,k}$ say that for large primes, the cohomology of $PW_{n,k}$ and $P_\ell W_{n,k}$ matches that of $Y_{n,k}$. For $W_{n,k;m}$, the expression matches the cohomology of a product of the lens space $L_m(2n-2k+1)$ and a bunch of odd dimensional spheres. Using elementary arguments, in this section we observe that these isomorphisms may be lifted to $p$-local homotopy equivalences for a sufficiently large lower bound on $p$.  

The first step towards a homotopical result starting from a cohomology isomorphism is a rational homotopy calculation. Note that the Serre spectral sequence for the fibration $W_{n,k}\longrightarrow PW_{n,k}\longrightarrow \C P^\infty  $ with rational coefficients tell us that 
\[
H^*(PW_{n,k};\Q)=\Q[x]/(x^{n-k+1})\otimes \Lambda_{\Q}(\gamma_{n-k+2},\cdots, \gamma_n), \]
 where $|x|=2, |\gamma_j|=2j-1 $.
\begin{mysubsection}{Rational splittings for $PW_{n,k}$} \label{ratdec}
The rational homotopy type of simply connected spaces are determined by its minimal model \cite{FHT01}. For a space $X$, we denote its minimal model by $m_X $. We shall show that $m_{PW_{n,k}}$ is isomorphic to  $m_{Y_{n,k}}$, which will tell us that they are rational homotopy equivalent.  

We know that 
$$m_{Y_{n,k}}= m_{\C P^{n-k}} \otimes m_{S^{2n-2k+3}} \otimes \cdots \otimes m_{S^{2n-1}}.$$
 Also $m_{S^{2j-1}}= (\Lambda (y_j),d=0)$, where $|y_j|=2j-1$, and $m_{\C P^{n-k}}= (P(x)\otimes \Lambda(y_{n-k+1}),d)$, where $|x|=2,|y_{n-k+1}|=2n-2k+1$ and $d(x)=0, d(y_{n-k+1})=x^{n-k+1}$ \cite{FHT01}. In this expression, $P(x)$ stands for the polynomial algebra on the generator $x$. 
\begin{prop}
The minimal model for $PW_{n,k}$ is given by 
\[
m_{PW_{n,k}}=P(\tilde{x})\otimes \Lambda(\tilde{y}_{n-k+1},\cdots,\tilde{y}_n),\]
 where $|\tilde{x}|=2, |\tilde{y}_j|=2j-1$ and the action of differential is determined by the following: $d(\tilde{x})=0, d(\tilde{y}_{n-k+1})=\tilde{x}^{n-k+1}, d(\tilde{y}_j)=0 $, $\forall n-k+1<j\leq n$.
\end{prop}
\begin{proof}
We only need to construct a differential graded algebra (DGA) morphism $\varphi$ from the the minimal Sullivan DGA  given in the statement to the DGA $A^*_{PL}(PW_{n,k})$ (following the same notation as \cite{FHT01})  which will also be a quasi-isomorphism. Actually it turns out that $PW_{n,k}$ is a formal space ie. its minimal model is determined by its cohomology ring which we will see below. Define $\varphi$ by the following action on generators
 \begin{myeq}\label{DGA-mor} \varphi(\tilde{x})= \omega, \varphi(\tilde{y}_{n-k+1})=\theta, \varphi (\tilde{y}_j)= \sigma_j,\forall n-k+1<j\leq n \end{myeq} where $\omega$ and $\sigma_j$ are cocycles in $A^*_{PL}(PW_{n,k})$ such that $[\omega]=x\in H^*(PW_{n,k};\Q) $, $[\sigma_j]=\gamma_j\in H^*(PW_{n,k};\Q)$. Since $x^{n-k+1}=0$ in $H^*(PW_{n,k};\Q)$, there will be some element $\theta \in A^{2n-2k+1}_{PL}(PW_{n,k})$ whose image under the differential in $A^*_{PL}(PW_{n,k})$ will be $\omega^{n-k+1}$. This is the $\theta$ used in \ref{DGA-mor}. With this definition, 
$\varphi$ is clearly a DGA quasi-isomorphism. 
\end{proof}
Now it is evident that both $PW_{n,k}$ and $Y_{n,k}$ have the same minimal model. Therefore, ${PW_{n,k}}_\Q \simeq {Y_{n,k}}_\Q$. 
\end{mysubsection}

\begin{mysubsection}{Splitting at large primes}
We now replicate the rational result of \S \ref{ratdec} in the $p$-local homotopy category. The rough estimate is so that $2p-3$ is larger than the dimension of the manifold, for which we obtain a result by elementary means. We note 
\[ \dim(PW_{n,k})= 2nk - k^2 - 1 = \dim(P_\ell W_{n,k}),~~ \dim(W_{n,k;m})= 2nk-k^2.\] 
In the following we use the class $x\in H^2(PW_{n,k})$ to obtain the map $PW_{n,k} \to \C P^\infty$. 
\begin{prop}\label{largpret}
For all primes $p> \frac{2nk-k^2-1}{2}+k-n$, there is a map $(PW_{n,k})_{(p)} \longrightarrow \C P^{n-k}_{(p)}$ which is a lift of the map  $(PW_{n,k})_{(p)} \to \C P^\infty_{(p)}$ up to homotopy. The same conclusion holds for $P_\ell W_{n,k}$ under the additional assumption that $p\nmid \sum_{|I|=n-k+1} l^I$.
\end{prop}
\begin{proof}
The homotopy fibre of the inclusion $\C P^{n-k} \to\C P^\infty$ is $S^{2n-2k+1}$. Therefore, obstructions for lifting the map $PW_{n,k}\to \C P^\infty$ on $p$-localizations lie in $H^r(PW_{n,k}; \pi_{r-1}(S^{2n-2k+1}_{(p)}))$, which are $0$ under the given hypothesis. As $\dim(PW_{n,k})< \dim(S^{2n-2k+1}) + 2p-3$, the only obstruction that may arise is if $r=2n-2k+2$, in which degree the cohomology of $PW_{n,k}$ is $0$ from \eqref{cohplarge}.   
%
%
The same argument also works for $P_\ell W_{n,k}$ under the additional hypothesis.
\end{proof}

Let $Z_{n,k} = \Big(S^{2n-2k+3} \times S^{2n-2k+5} \times \cdots \times S^{2n-1}\Big)_{(p)}$. The cohomology of $Z_{n,k}$ is the exterior algebra
\[H^\ast(Z_{n,k}) \cong \Lambda(\epsilon_{2n-2k+3},\epsilon_{2n-2k+5}, \cdots, \epsilon_{2n-1}),\]
where the classes $\epsilon_{2j-1}$ are pullbacks of the generators of $H^\ast(S^{2j-1})$ via the corresponding projection.   We prove  
\begin{prop}\label{largepsph}
Suppose $p$ is as in the hypothesis of Proposition \ref{largpret}. There is a map $\rho:(PW_{n,k})_{(p)} \to Z_{n,k} $ such that $\rho^\ast(\epsilon_{2j-1})=\gamma_j $ for all $n-k+2\leq j \leq n$, with $\gamma_j$ as in \eqref{choicegen}. The same result holds for $P_\ell W_{n,k}$ if $p\nmid \sum_{|I|=n-k+1} l^I$.
\end{prop}
\begin{proof}
Suppose  $L= K(\Z_{(p)},2n-2k+3)\times \cdots \times K(\Z_{(p)},2n-1) $, then we get a map $\nu: Z_{n,k} \to L $ which classifies the cohomology generators of $Z_{n,k}$. Now, 
\[\pi_r(Z_{n,k}) \cong\pi_r(S^{2n-2k+3}_{(p)}) \oplus \cdots \oplus \pi_r(S^{2n-1}_{(p)}),\]  
 and $\pi_r(S^{2n-2k+j}_{(p)})=0$ for $2n-2k+j<r<2n-2k+j+2p-3$, imply that $\pi_r(Z_{n,k})= \Z_{(p)}$, for $r$ odd and $r= 2n-2k+3, \cdots, 2n-1$ and $\pi_r(Z_{n,k})=0$, for all other $r\leq \dim(PW_{n,k})$. So $\nu$ must be a $\dim(PW_{n,k})$-equivalence. 

From the conclusion of the paragraph above, we will get the following isomorphism 
$$ [(PW_{n,k})_{(p)}, Z_{n,k}] \stackrel{\nu_\ast}{\to} [(PW_{n,k})_{(p)},L]\cong  H^{2n-2k+3}(PW_{n,k};\Z_{(p)})\oplus \cdots \oplus H^{2n-1}(PW_{n,k};\Z_{(p)}) .$$ 
Hence we get a map $\rho: (PW_{n,k})_{(p)}\to Z_{n,k}$ that pulls back the cohomology generators  $\epsilon_{2j-1}\in H^\ast(Z_{n,k})$ corresponding to the generators of $H^*(S^{2j-1}_{(p)})$ to $\gamma_j\in H^\ast(PW_{n,k},\Z_{(p)})$. This argument works entirely analogously for $P_\ell W_{n,k}$.
\end{proof}

We may now assemble the two results from Propositions \ref{largpret} and \ref{largepsph} to get  maps for $p$ large 
\[ 
\Big(PW_{n,k}\Big)_{(p)} \to \Big(Y_{n,k}\Big)_{(p)}, ~~\Big(P_\ell W_{n,k}\Big)_{(p)} \to \Big(Y_{n,k}\Big)_{(p)}
\]
which are cohomology isomorphisms. Thus, we have proved the following result. 
\begin{thm} \label{largepdec}
Suppose $p>  \frac{2nk-k^2-1}{2}+k-n$. Then we have, 
\[ \Big(P W_{n,k}\Big)_{(p)} \simeq \Big[ \C P^{n-k} \times S^{2n-2k+3} \times \cdots \times S^{2n-1} \Big]_{(p)}.\]
If further $p\nmid \sum_{|I|=n-k+1} l^I$,
\[ \Big(P_\ell W_{n,k}\Big)_{(p)} \simeq \Big[ \C P^{n-k} \times S^{2n-2k+3} \times \cdots \times S^{2n-1} \Big]_{(p)}.\]
\end{thm}  

Theorem \ref{largepdec} may be used to provide an $S^1$-equivariant decomposition of the Stiefel manifold. 
\begin{prop} \label{eqsplitlarg}
For $p>\frac{2nk-k^2-1}{2}+k-n$, we have the following splitting as $S^1$-spaces. 
$$ \Big(W_{n,k}\Big)_{(p)} \simeq [S^{2n-2k+1} \times S^{2n-2k+3}\times \cdots \times S^{2n-1}]_{(p)}.$$ 
\end{prop}
\begin{proof} 
For the bound on $p$ stated in the proposition, there is a map $\phi_k:(PW_{n,k})_{(p)} \to \C P^{n-k}_{(p)}$ lifting the classifying map of the $S^1$-bundle $W_{n,k}\to PW_{n,k}$. Hence we have the following homotopy pullback diagram, 
$$ \xymatrix{(W_{n,k})_{(p)}\ar[r]\ar[d] & S^{2n-2k+1}_{(p)}\ar[d]\\
 (PW_{n,k})_{(p)}\ar[r] & \C P^{n-k}_{(p)} } $$
leading to an $S^1$-equivariant map $(W_{n,k})_{(p)}\to S^{2n-2k+1}_{(p)}$. 
We also have an $S^1$-equivariant map $(W_{n,k})_{(p)} \to (W_{n,k-1})_{(p)}$ arising from the $S^1$-equivariant projection $W_{n,k} \to W_{n,k-1}$. The product of these two maps is evidently a cohomology isomorphism. So we obtain an $S^1$-equivariant equivalence 
$$ (W_{n,k})_{(p)} \simeq S^{2n-2k+1}_{(p)} \times (W_{n,k-1})_{(p)}.$$
 And by induction on $k$, we arrive at the result stated in the proposition.
\end{proof}
Now observe that both $W_{n,k;m}$ and $L_m(2n-2k+1)\times S^{2n-2k+3}\times \cdots \times S^{2n-1}$ have the same $\Z_{(p)}$-cohomology whenever $p>n$. 
\begin{prop}\label{wnklarg}
For $p>\frac{2nk-k^2-1}{2}+k-n$ we have the the following splitting 
$$ \Big(W_{n,k;m}\Big)_{(p)} \simeq  [L_m(2n-2k+1)\times S^{2n-2k+3}\times \cdots \times S^{2n-1}]_{(p)}.$$
\end{prop}
\begin{proof}
For $p>\frac{2nk-k^2-1}{2}+k-n$, we have the $S^1$-equivariant map $W_{n,k} \to S^{2n-2k+1}$ as shown in Proposition \ref{eqsplitlarg}. So we get a map $W_{n,k;m} \to L_m(2n-2k+1)$ by considering the $C_m$-orbit spaces. Now comparing the spectral sequences associated to the fibrations 
$$S^1 \to W_{n,k;m} \to PW_{n,k}, \mbox{ and } S^1 \to L_m(2n-2k+1) \to \C P^{n-k},$$ 
we can see that the map $W_{n,k;m} \to L_m(2n-2k+1)$ induces an isomorphism on $ H^j(- ;\Z_{(p)}) $ for $j \leq 2n-2k+1$. We also have the maps $W_{n,k;m} \to PW_{n,k} \to S^{2n-2k+2r+1}$ for $ k>r>0$ and they map the cohomology generators of $S^{2j-1}$ to the corresponding generator of $H^\ast(W_{n,k;m})$ in \eqref{wnkmcoh}. Hence we have a map from $W_{n,k;m}$ to $ L_m(2n-2k+1)\times S^{2n-2k+3}\times \cdots \times S^{2n-1} $ that induces an isomorphism on cohomology, and both spaces being simple, we get the desired equivalence. 
\end{proof}
\end{mysubsection}

\section{Stable decompositions of projective Stiefel manifolds} \label{decst}

In this section, we prove stable decomposition results for the projective Stiefel manifold $PW_{n,k}$ at primes greater than $n$. The crucial observation here is that if $p > n$, we have the stable equivalence 
\[ 
\C P^n_{(p)} \simeq S^2_{(p)} \vee S^4_{(p)} \vee \cdots \vee S^{2n}_{(p)}. 
\] 
This fact may be proved by showing that the attaching maps of the cells of the $p$-local complex projective space are stably trivial if $p>n$. The reason is that the first non-trivial $p$-torsion in the stable homotopy groups of $S^0$ occurs in degree $2p-3$, which is greater than the degrees of the attaching maps of $\C P^n$ if $p>n$. However one requires a more delicate argument to make this work for $PW_{n,k}$, because there are cell attachments of degree greater than $2p-3$. 
 
\begin{mysubsection}{Minimal cell structures} \label{mincell}
We recall the minimal cell structures for CW complexes from \cite[\S 4.C]{Hat01}. For a simply connected CW complex with finitely generated homology groups, we have a CW complex structure with ``minimum number of cells". More precisely, writing the homology groups using generators and relations, we have one ``generator" $n$-cell for every generator of a cyclic copy of $H_n(X)$, and one ``relator" $n+1$-cell whose boundary is $k$ times a ``generator" $n$-cell in cellular homology whenever the cyclic copy in $H_n(X)$ corresponding to the generator is $\Z/k$.    

We will use a version of the above result for $p$-local spaces. In this case we use $p$-local cells $\DD^n_{(p)}$ which is the cone of $S^{n-1}_{(p)}$, and a $p$-local $n$-cell attachment to a $p$-local space $X$ is the attachment of $\DD^n_{(p)}$ along a map $S^{n-1}_{(p)}\to X$. A $p$-local finite CW complex is one which is obtained by attaching finitely many $p$-local cells in increasing dimensions. The following proposition follows directly using the arguments of \cite[Proposition 4.C.1]{Hat01}. 
\begin{prop}\label{p-min-cell}
Suppose that $X$ is a simply connected $p$-local space such that $H_\ast (X;\Z_{(p)})$ is finitely generated and torsion free. Then, $X$ has a $p$-local CW complex structure with one $p$-local $n$-cell for each basis element of $H_n(X;\Z_{(p)})\cong \Z_{(p)}^l$.
\end{prop}

\begin{rmk}\label{min-cell-att}
Suppose that $X$ is as in Proposition \ref{p-min-cell}. We observe that as the boundary map in  cellular homology is $0$, the attaching map of such a cell factors through the $(n-2)$-skeleton. 
\end{rmk}
\end{mysubsection}

\begin{mysubsection}{The Chern character} \label{chchar}
Recall the Chern character \cite[Chapter 5]{Hil71} 
\[
ch : K^\ast(X) \to H^\ast(X;\Q)[u^{\pm}]
\]
where $u$ is a degree $2$ class inducing the $2$-periodicity in the second factor. We restrict our attention to $p$-local finite CW complexes of Proposition \ref{p-min-cell}, and note that the arguments in page 73 of \cite{Hil71} implies the following result. 
\begin{prop}\label{p-ch-char}
Suppose that $X$ is a simply connected $p$-local space such that $H_\ast (X;\Z_{(p)})$ is finitely generated and torsion free. Then, $K^\ast(X)\otimes \Z_{(p)}$ is torsion-free and the Chern character induces an isomorphism $K^\ast(X)\otimes \Q \cong H^\ast(X;\Q)[u^{\pm}]$. 
\end{prop}
\end{mysubsection}

\begin{mysubsection}{A criterion for stable splitting} \label{stsplitcrit}
We now formulate a general criterion to have a $p$-local stable decomposition into a wedge of spheres. This involves ruling out attaching maps in the image of the $J$-homomorphism using an assumption on the Chern character. For the remainder of this section, we work in the category $\Sa$ which stands for the stable homotopy category, whose objects are sequential spectra and morphisms are the homotopy classes of maps between them. For a space $X$, we use the same notation for the suspension spectrum as an object of $\Sa$. We also use the Chern character for cellular spectra with finitely many cells. Note from \cite[Theorem 1.1.14]{Rav86} that elements in $\pi_k^s(S^0)$ for $k>0$ lie in the image of the $J$-homomorphism if $k<2p^2-2p$.   
\begin{theorem} \label{resstsplit}
Let $X$ be a  simply connected $p$-local  finite CW-complex satisfying the following conditions 
\begin{enumerate}
    \item $\dim X <2p^2-2p$.
    \item $H^*(X;\Z_{(p)})$ is free as a $\Z_{(p)} $-module.
    \item The Chern character map for $X$ has image in  $H^*(X;\Z_{(p)}) \subset H^*(X;\Q )$.
\end{enumerate}
Then, $X \simeq$ a wedge of $p$-local spheres. 
\end{theorem}
\begin{proof}
We consider the minimal cell structure of $X$ of Proposition \ref{p-min-cell} and prove that the attaching maps are $0$ in the stable homotopy category. We show this by induction over $k$ for $X^{(k)}$, which is the $k^{th}$-skeleton of $X$. The properties of the minimal cell structure, the naturality of the Chern character together with Proposition \ref{p-ch-char}, implies that any sub-complex $A$ of $X$ satisfies the three hypotheses stated in the theorem. 

We assume that $X^{(k)}$ splits as a wedge of $p$-local cells, and we want to conclude the same for $X^{(k+1)}$. Since there are finitely many cells, it suffices to prove that every attaching map $\phi: S^k_{(p)} \longrightarrow X^{(k)}$ is trivial. We write $Y$ for the mapping cone on $\phi$. From Remark \ref{min-cell-att}, this attaching map will actually land in $X^{(k-1)}_{(p)}$. 
Writing
$$X^{(k)} \simeq S^{n_1}_{(p)} \vee \cdots \vee S^{n_r}_{(p)},$$
 the attaching map of the $p$-local $(k+1)$-cell in $X_{(p)}$ must factor through $\bigvee_{n_j \neq k} S^{n_j}_{(p)}$, and is therefore a sum of maps of the form $f_j: S^k_{(p)}\to S^{n_j}_{(p)}$. Note that $k-n_j<2p^2-2p$ by condition 1), so that  the homotopy class of these maps  lie in the image of $J$ homomorphism. Since these classes are in odd degree, we may assume $k-n_j$ is odd. 
Now we have the following homotopy commutative diagram
\begin{myeq} \label{proj}
 \xymatrix{S^k_{(p)}\ar[r]^{\phi}\ar[d] &X^{(k)}_{(p)}\ar[r]\ar[d] & Y_{(p)}\ar[r]\ar[d] &S^{k+1}_{(p)}\ar[d]^{=}\\ 
S^k_{(p)}\ar[r]_-{f_j} &S^{n_j}_{(p)}\ar[r] & C(f_j)\ar[r] &S^{k+1}_{(p)} } 
\end{myeq}
where the second vertical map from left is induced by the  retraction
$$S^{n_1}_{(p)} \vee \cdots \vee S^{n_r}_{(p)} \longrightarrow S^{n_j}_{(p)}. $$
To show that $f_j$ is trivial it suffices to show that the $K$-theoretic $e$-invariant \cite[\S 7]{Ada66} of $f_j$, $e_K(f_j)$ vanishes. Note that \cite[Proposition 7.14]{Ada66} implies that we may compute this using complex $K$-theory as $p\neq 2$. We also notice that since we are dealing with $p$-local spheres, the $e$-invariant will take values in $\Q/\Z_{(p)}$. 

We look at the following diagram of short exact sequences induced by the Chern character:
$$ \xymatrix{0&\ar[l]\ar[d]_-{ch} K^*_{(p)}(X^{(k)})& \ar[l]\ar[d]_-{ch} K^*_{(p)}( Y) &\ar[l]\ar[d]^-{ch} K^*_{(p)}(S^{k+1}) &\ar[l] 0\\
0 &\ar[l]H^{*}(X^{(k)};\Q)[u^\pm] &\ar[l] H^{*}( Y;\Q)[u^\pm] &\ar[l] H^{*}(S^{k+1};\Q)[u^\pm] &\ar[l]0  } $$
where the short exactness of the first row follows from the short exactness of the bottom row along with the injectivity of the Chern character map for $p$-local sphere. As $Y$ is a subcomplex of $X$, the image of the Chern character lies in $H^\ast(Y; \Z_{(p)})$ by condition 3). Now from the following diagram (the injectivity of the horizontal arrows follow from \eqref{proj})
\[
\xymatrix{K_{(p)}(C(f_j))\ar@{^{(}->}[r] \ar[d]_-{ch} &K_{(p)}( Y)\ar[d]^-{ch}\\ 
H^*(C(f_j);\Q)[u^\pm]\ar@{^{(}->}[r] &H^*(Y;\Q)[u^\pm], }
\]
we conclude that the image of the Chern character for $C(f_j)$ also lies in $\Z_{(p)}$-cohomology. Applying the reformulation of the $e$-invariant in \cite[Proposition 7.8]{Ada66}, we deduce that $e(f_j)=0$, and hence $f_j\simeq 0$. 
%
%
%
\end{proof}
\end{mysubsection}

\begin{mysubsection}{Algebra generators for $K^*(W_{n,k})$ and $K^\ast(PW_{n,k})$}
We now write down the structure of $K^\ast PW_{n,k}$ (and $K^\ast W_{n,k}$), identifying algebra generators as was done in \eqref{choicegen} for ordinary cohomology. The coefficient ring for $p$-local complex $K$-theory is $\Z_{(p)}[\beta^\pm]$ where $\beta$ stands for the Bott element lying in degree $2$. Recall that complex $K$-theory is complex oriented, and $x_K=\beta^{-1} (L-1)$ is a choice of complex orientation for the canonical line bundle $L$ over $\C P^\infty$. This implies that every complex vector bundle is $K$-orientable, and defines universal Chern classes $c_j^K$ such that $K_{(p)}^\ast(BU(n)) \cong \Z_{(p)}[\beta^\pm][[c_1^K,\cdots, c_n^K]]$.  

We note that the complex orientation expressed as a map 
\[ x_K: \C P^\infty \to \Sigma^2 K_{(p)}\]
may be expressed as a composite via the connective cover $ku_{(p)}$ of $K_{(p)}$
\[ \C P^\infty  \stackrel{x_{ku}}{\to} \Sigma^2 ku_{(p)} \to \Sigma^2 K_{(p)}\]
as $\C P^\infty$ is $1$-connected. This implies that the universal Chern classes $c_j^K$ are the image of $ku_{(p)}$-Chern classes $c_j^{ku}$. 
Following the method described in \cite[Proposition 3.4 and Proposition 3.6]{BD21} we obtain classes $\tau_j^{ku} \in ku_{(p)}^*(BU(n))$ such that $\tau_j^{ku}=c_j^{ku}+\sum_{k>j}\nu_k c_k^{ku}$, where $\nu_k \in ku_{(p)}^*(pt)$. We use the map $\alpha: (BU(n-k),W_{n,k})\to (BU(n),\ast)$ arising from the fibration $W_{n,k}\to BU(n-k)\to BU(n)$, and the diagram  
$$\xymatrix{\cdots \ar[r] &ku_{(p)}^i(W_{n,k}) \ar[r]^-\delta &ku_{(p)}^{i+1}(BU(n-k),W_{n,k})\ar[r]  &ku_{(p)}^{i+1}(BU(n-k))\ar[r]  &\cdots \\
&&\widetilde{ku_{(p)}}^{i+1}(BU(n)).\ar[u]^-{\alpha^\ast}}$$ 
to define $y_j^{ku}$ by the equation $\delta (y_j^{ku})=\alpha^\ast (\tau_j^{ku})$ for $j\geq n-k+1$. The classes $y_j^{ku}$ serve as generators of $ku_{(p)}^*(W_{n,k})$, that is, 
\[ ku_{(p)}^\ast (W_{n,k})\cong \Lambda_{ku_{(p)}^\ast} (y_{n-k+1}^{ku},\cdots, y_n^{ku}).\]
Pushing forward to $K_{(p)}^*$, by the map $ku_{(p)} \to K_{(p)}$, we obtain classes $\tau_j^K \in K_{(p)}^*(BU(n))$ such that $\tau_j^K=c_j^K+\sum_{k>j}\nu_k c_k^K$, and $y_j^K \in K_{(p)}^*(W_{n,k})$ such that in the diagram 
$$\xymatrix{\cdots \ar[r] &K_{(p)}^i(W_{n,k}) \ar[r]^-\delta &K_{(p)}^{i+1}(BU(n-k),W_{n,k})\ar[r]  &K_{(p)}^{i+1}(BU(n-k))\ar[r]  &\cdots \\
&&\widetilde{K_{(p)}}^{i+1}(BU(n)),\ar[u]^-{\alpha^\ast}}$$  
$\delta(y_j^K)= \alpha^\ast (\tau_j^K)$, and one can show $y_j^K$ give a set of generators for $K_{(p)}^*(W_{n,k})$. We now multiply $y_j^K$ by a power of $\beta$ so that  $y_j^K \in K_{(p)}^{-1}(W_{n,k})$. For the projective Stiefel manifold we have 
$$K_{(p)}^*(PW_{n,k})\cong K_{(p)}^\ast \otimes_{BP^\ast} BP^\ast(PW_{n,k}) \cong \Lambda_{K_{(p)}^*(pt)}(\gamma_{n-k+3},\cdots, \gamma_n)\otimes _{K_{(p)}^*(pt)} K_{(p)}^*(pt)[[x]]/I$$
 where $|\gamma_i|=-1$ and $|x|=0$ and $I$ is the ideal generated by $\{\binom{n}{j}x^j|n-k <j \leq n \}$ \cite{BD21}.  We are assuming $p>n$,  so that 
$$K_{(p)}^*(PW_{n,k})\cong \Lambda_{K_{(p)}^*(pt)}(\gamma_{n-k+3},\cdots, \gamma_n)\otimes _{K_{(p)}^*(pt)} K_{(p)}^*(pt)[[x]]/(x^{n-k+1}) $$ where $|\gamma_i|=-1$ and $|x|=0 $. We may construct the $K$-theoretic algebra generators $\gamma_j^K$ via the following analogue of \eqref{PW-H}
\begin{myeq}\label{PWK}
\xymatrix{0\ar[d]\\K_{(p)}^{-1}(PW_{n,k})\ar[d]_-{\delta} &0\ar[d]\\K_{(p)}^0 (\C P^{\infty},PW_{n,k})\ar[d]_-{i^*} &K_{(p)}^0 (BU(n),BU(n-k))\ar[l]^-{f^*}\ar[d]\\
K_{(p)}^0 (\C P^{\infty}) &K_{(p)}^0(BU(n))\ar[l]^-{f_0^*}}
\end{myeq}
Identify $K_{(p)}^0(BU(n),BU(n-k))$ with the ideal in $K_{(p)}^0(BU(n)) $ generated by $c_j^K$ for $n-k<j\leq n$, and write $u_j^K=f^*\tau_j^K$. Since $$f_0^*\tau_j^K=\binom{n}{j}x_K^j + \sum_{k>j} \binom{n}{k}\nu_k x_K^k= x_K^j \mu_j^K,$$
 where $\mu_j^K$ is a unit in $K_{(p)}^\ast (\C P^{\infty})$, we see that 
$$i^*(u_j^K - x_K^{j-(n-k+1)}\frac{\mu_j^K }{\mu_{n-k+1}^K} u_{n-k+1}^K)=0.$$ 
Let $\gamma_j^K \in K_{(p)}^{-1}(PW_{n,k})$ be defined by 
\begin{myeq} \label{choicegenK}
 \delta \gamma_j^K = u_j^K - x_K^{j-(n-k+1)}\frac{\mu_j^K }{\mu_{n-k+1}^K} u_{n-k+1}^K.
\end{myeq}
In the following proposition, we show that $\gamma_j^K$ generates the exterior algebra part of $K^*_{(p)}(PW_{n,k})$. 
\begin{prop}\label{genkprojst}
With notations as above, 
$$K_{(p)}^*(PW_{n,k})\cong \Lambda_{K_{(p)}^*(pt)}(\gamma_{n-k+3}^K,\cdots, \gamma_n^K)\otimes _{K_{(p)}^*(pt)} K_{(p)}^*(pt)[x]/(x^{n-k+1}) . $$
\end{prop}
\begin{proof}
Consider the homotopy fixed point spectral sequence \cite[Proposition 2.1]{BD21}
\[ E_2^{s,t} = \Z[y]\otimes K_{(p)}^t(W_{n,k}) \implies K_{(p)}^{s+t}(PW_{n,k}),\]
and from the analogous computation of \cite[Proposition 4.5]{BD21} deduce that a set of elements form generators of the exterior algebra part if they pullback to corresponding exterior algebra generators of $K_{(p)}^\ast(W_{n,k})$. Consequently, it suffices to prove that $\pi^*\gamma_j^K=y_j^K$ in $K^*_{(p)}(W_{n,k})$. From the diagram \eqref{pullback} we get the following commutative diagram extending \eqref{PWK}
$$ \xymatrix{K^{-1}_{(p)}(W_{n,k})\ar[d]_-{\delta}^-{\cong}&K^{-1}_{(p)}(PW_{n,k})\ar[d]^-{\delta}\ar[l]_-{\pi^*}\\
                         K_{(p)}^0(E,W_{n,k}) &K_{(p)}^0(\C P^{\infty},PW_{n,k})\ar[l]^-{\pi^*} &K_{(p)}^0(BU(n),BU(n-k)).\ar[l]^-{f^*}} $$
 Here we compute for $j>n-k+1$,
 \begin{align*}
 \delta\pi^*(\gamma_j^K) &= \pi^*\delta(\gamma_j^K) \\
                                            &= \pi^*(u_j^K - x_K^{j-(n-k+1)}\frac{\mu_j^K }{\mu_{n-k+1}^K} u_{n-k+1}^K) \\
                                            &= \pi^*u_j^K \mbox{ (as } \pi^\ast(x)=0)\\
                                            & = \pi^* f^*(\tau_j^K). 
 \end{align*}  
 
 On the other hand, we have the following diagram
  $$ \xymatrix{&U(n)\ar[rr]\ar[ld]\ar|(0.5){\hole}[ldd]  &&\ast \ar[ld]\ar[ldd]\\
                           W_{n,k}\ar[rr]\ar[d] &&BU(n-k)\ar[d]\\
                           E\ar[rr] &&BU(n)}$$
 which leads to the diagram below
 $$ \xymatrix{K_{(p)}^{-1}(W_{n,k})\ar[r]^-{\cong}\ar@{^{(}->}[d] &K_{(p)}^0(E,W_{n,k})\ar[d] &K_{(p)}^0(BU(n),BU(n-k))\ar[l]\ar[d]\\
                          K_{(p)}^{-1}(U(n))\ar[r]_-{\cong} &K_{(p)}^0(E,U(n)) &K_{(p)}^0(BU(n),\ast)\ar[l]} $$
In the bottom row of above diagram image of $y_j^K$ and image of $\tau_j^K$ coincide by the construction of $y_j^K$. So in the first row the same will happen for $j> n-k$. Hence we must have $\delta \pi^* \gamma_j^K = \delta y_j^K$, and $\delta$ being injective, we get our desired result. 
\end{proof}
\end{mysubsection}

\begin{mysubsection}{Stable decompositions for ${PW_{n,k}}_{(p)}$ } 
We verify the hypothesis of Theorem \ref{resstsplit} for $PW_{n,k}$ if $n<p$. First, we prove a lemma regarding the image of the Chern character. Recall that 
$$H^*(PW_{n,k};\Z_{(p)})\cong \Lambda_{\Z}(\gamma^H_{n-k+2},\cdots,\gamma^H_{n})\otimes \Z_{(p)}[x]/(x^{n-k+1})$$
 with $|\gamma^H_j|=2j-1$ as in \eqref{choicegen}, $|x|=2$. 
In the following lemma, we use the construction of $\gamma_j^K$ of \eqref{choicegenK} which generate the exterior algebra part of the $K$-theory of ${PW_{n,k}}_{(p)}$ by Proposition \ref{genkprojst}. 
\begin{lemma}\label{chgammaexp}
For $n-k+2 \leq j \leq n$ , $ch(\gamma_j^K)$ is a $\Q$-linear combination of $\{\gamma^H_s x^t \mid n-k+2\leq s \leq n , ~0\leq t \leq n-k  \}$. 
\end{lemma}
\begin{proof}
Consider the commutative diagram
$$\xymatrix{K^{-1}(PW_{n,k})\ar@{^{(}->}[r]^{\delta} \ar[d]^-{ch} &K(\C P^{\infty},PW_{n,k})\ar[d]^-{ch} &K(BU(n),BU(n-k))\ar[l]^{(f_0,f)^\ast}\ar[d]^{ch}\\
H^{odd}(PW_{n,k};\Q) \ar@{^{(}->}[r]^-{\delta} &H^{ev}(\C P^{\infty},PW_{n,k};\Q) &H^{ev}(BU(n),BU(n-k);\Q)\ar[l]^-{(f_0,f)^\ast}}  $$
with $f$, $f_0$ as defined in \eqref{pullback}. Now $H^\ast(BU(n),BU(n-k);\Q)$ is the ideal $(c_{n-k+1},\cdots,c_n)$ of $H^\ast(BU(n);\Q)=\Q[c_1,\cdots,c_n]$. In these terms, we must have an equation of the form
$$ch(\tau_j^K)=c_{n-k+1}\cdot P_{n-k+1,j}+c_{n-k+2}\cdot P_{n-k+2,j} \cdots+ c_n \cdot P_{n,j}$$ 
for some power series $P_{i,j} \in \Q[[c_1,\cdots,c_n]] $.  The right commutative square in the above implies (as $(f_0,f)^\ast$ preserves the multiplication of relative cohomology classes)
$$ ch(u_j^K)= ch(f^\ast \tau_j^K) = f^*(c_{n-k+1})\cdot f_0^*(P_{n-k+1,j})+f^\ast(c_{n-k+2})\cdot f_0^\ast (P_{n-k+2,j}) +\cdots + f^*(c_n)\cdot f_0^*(P_{n,j})$$ 
with $f_0^*(P_{i,j})\in H^\ast(\C P^{\infty};\Q)=\Q[[x]]$. We may now compute from \eqref{choicegenK}
\begin{align*}
\delta ch(\gamma_j^K) & = ch(\delta \gamma_j^K) \\ 
                                  &= ch(u_j^K - x_K^{j-(n-k+1)}\frac{\mu_j^K }{\mu_{n-k+1}^K}  u_{n-k+1}^K) \\ 
                                 &= ch(f^\ast \tau_j^K - x_K^{j-(n-k+1)}\frac{\mu_j^K }{\mu_{n-k+1}^K} f^\ast \tau_{n-k+1}^K) \\
                                 &= f^*(c_{n-k+1})\cdot f_0^*(Q_{n-k+1,j})+\rho_{n-k+2}\cdot f_0^*(Q_{n,j}) +\cdots + \rho_n\cdot f_0^*(Q_{n,j}),
\end{align*}
where $\rho_j$ is defined in \eqref{choicegen} as, 
\[ 
\rho_j=  u_j^H - x^{j-(n-k+1)}\frac{\mu_j}{ \mu_{n-k+1}} u_{n-k+1}^H = f^\ast(c_j) - x^{j-(n-k+1)}\frac{\mu_j}{ \mu_{n-k+1}} f^\ast(c_{n-k+1}),
\]
and for some polynomials $Q_{i,j}$ in the $c_i$. The element $\delta ch(\gamma_j^K)$ maps to $0$ in $H^\ast (\C P^\infty)$ from which it follows that $f_0^\ast Q_{n-k+1,j}$ must be $0$, so that 
\[
\delta ch(\gamma_j^K) = f^*(\rho_{N+1})\cdot f_0^*(Q_{N+1,j}) +\cdots + f^*(\rho_n)\cdot f_0^*(Q_{n,j}). 
\]
We now write $f_0^\ast(Q_{i,j})= \phi_{i,j}(x)$, and use the fact that $\delta : H^\ast(PW_{n,k}) \to H^{\ast +1 }(\C P^\infty, PW_{n,k})$ is a map of $H^\ast(\C P^\infty)$-modules, which implies 
\begin{align*}
\delta ch(\gamma_j^K) & =\rho_{n-k+2}\cdot f_0^*(Q_{N+1,j}) +\cdots + \rho_n\cdot f_0^*(Q_{n,j})\\
                                   & =  \delta(\gamma_{n-k+2}^H)\cdot \phi_{N+1,j}(x) +\cdots + \delta(\gamma_n)\cdot \phi_{n,j}(x)\\
                                   &= \delta(\gamma_{n-k+2}^H \phi_{N+1,j}(x) +\cdots + \gamma_n^H\cdot \phi_{n,j}(x)). 
\end{align*}
As $\delta$ is injective, we have that $ch(\gamma_j)$ is a linear combination of $\gamma_j^H x^r$ where $n-k+2\leq j \leq n$ and $r\leq n-1$. Further as $x^{n-k+1} $ maps to $0$ in $H^*(PW_{n,k};\Q)$, $ch(\gamma_j^K)$ must be a $\Q$-linear combination of the set described in the statement of the lemma.
\end{proof}

We now have all the ingredients in place to prove that the projective Stiefel manifold splits into a wedge of spheres in the stable homotopy category. 
\begin{thm}\label{projstsplit}
Let $p$ be a prime $> n$. Then, the $p$-localization of the projective Stiefel manifold ${PW_{n,k}}_{(p)}$ stably splits as a wedge of $p$-local spheres.   
\end{thm}
\begin{proof}
We verify that $PW_{n,k}$ satisfies the conditions of Theorem \ref{resstsplit}. The condition (2) already follows from the expression of \eqref{cohplarge}. We also have
\begin{align*}
\dim(PW_{n,k}) &= 2nk-k^2 - 1 = 2n(k-\frac{k^2}{2n}) -1 \\ 
                        & < 2n(n-1) \mbox{ as } k \leq n-1 \\ 
                        & < 2p(p-1) \mbox{ as } n<p,                        
\end{align*} 
which verifies the condition (1). Finally, we have to verify that the Chern character $ch$ has image in $H^\ast(PW_{n,k};\Z_{(p)})$. As $ch$ is a map of rings, it suffices to verify this on the generators $x$ and $\gamma_j^K$ of Proposition \ref{genkprojst}. The class $x$ is the pullback of the complex orientation class via $PW_{n,k}\to \C P^\infty$, so we have 
\[
ch(x) = e^x -1 = \sum_{i=1}^{n-k} \frac{x^i}{i!} \mbox{ as } x^{n-k+1}=0,
\]
which clearly lies in $H^\ast(PW_{n,k};\Z_{(p)})$ as $p>n-k$. The Chern character on the classes $\gamma_j^K$ is described by Lemma \ref{chgammaexp} on which we now apply the integrality result of Adams \cite[Theorem 1]{Ada61}. As $W_{n,k}$ is $2(n-k)$-connected, from the fibration 
\[ 
W_{n,k}\to PW_{n,k}\to \C P^{\infty}
\]
we obtain that there is a cell structure on $PW_{n,k}$ whose $2(n-k)$-skeleton is homotopy equivalent to $\C P^{n-k}$. It follows that the restriction of the classes $\gamma_j^K$ to $PW_{n,k}^{(2(n-k))}$ is $0$ with respect to this cell structure. The expression of Lemma \ref{chgammaexp} implies that the highest degree term which may occur in the expression of $ch(\gamma_j^K)$ is $\gamma_n^H x^{n-k}$, and this lies in degree $2n-1+2(n-k)$. Therefore, following \cite{Ada61}, we have  $ch_{2(n-k)+1+2r}(\gamma_j^K)=0$ if $r\geq n-1$. For $r<n-1<(p-1)$, $m(r)$ is not divisible by $p$, so by \cite[Theorem 1]{Ada61}, the proof is complete.  
\end{proof}
%

\begin{rmk}\label{splplwnk}
The argument above may be easily modified to deduce that the spaces $P_\ell W_{n,k}$ have a $p$-local stable decomposition into a wedge of spheres under the condition $p>n$, and $\sum_{|I|=n-k+1} l^I$ is not divisible by $p$. The latter condition implies that the cohomology of $P_\ell W_{n,k}$ with $\Z_{(p)}$-coefficients is torsion-free by \eqref{plwnkcoh}. The $K$-theory is also torsion-free, and the generators may be chosen via the diagram \eqref{pullbackpl} using the pair $(BU(k),Gr_k(\C^n))$ in place of $(BU(n), BU(n-k))$ in the calculations above. An analogous calculation implies that the Chern character has image in $\Z_{(p)}$ and then, we realize the $p$-local decomposition using Theorem \ref{resstsplit}.   
\end{rmk}

\end{mysubsection}

\section{Unstable $p$-local decompositions}\label{unpld}
In this section, we use the stable decomposition of the $p$-local $PW_{n,k}$ of \S \ref{decst} for $p>n$, to deduce an unstable product decomposition. Throughout this section, we work in the category of $p$-local spaces. As in \S \ref{decvlarge}, we compare $PW_{n,k}$ with the product 
 $Y_{n,k}=\C P^{n-k}\times S^{2(n-k)+3}\times S^{2(n-k)+5}\times \cdots \times S^{2n-1}$, which also stably decomposes into a wedge of spheres  for $p>n$.  Note that if $p>n$,  $PW_{n,k}$ and $Y$ both have the same cohomology with $\Z_{(p)}$-coefficients. 

In order to show the equivalence between $PW_{n,k}$ and $Y_{n,k}$, we construct maps from $PW_{n,k}$ to the individual factors of $Y_{n,k}$ in the $p$-local category. This is possible for small values of $k$ in comparison to $p$ and $n$. We first construct maps from $PW_{n,k}$ to $\C P^{n-k}$ which lift the map $PW_{n,k} \to \C P^\infty\simeq K(\Z,2)$ classifying the class $x\in H^2(PW_{n,k})$. 
\begin{prop}\label{retcp}
Suppose that $p>n+1$ and $k \leq \text{min}\{n,p+n-\sqrt{p^2+n^2-2p+1} \}$. Then, the map   $PW_{n,k}\to \C P^\infty$  lifts to $\C P^{n-k}$.
\end{prop}
\begin{proof} 
For $n=k$, there is nothing to prove. 
 We assume $n>k$, and that we have a $p$-local minimal cell structure on $PW_{n,k}$ via Proposition \ref{p-min-cell}. The $\Big(2(n-k)+2\Big)$-skeleton of $PW_{n,k}$ is homotopy equivalent to $\C P^{n-k}$. It suffices to prove that this equivalence extends all the way to a map $PW_{n,k}\to \C P^{n-k}$. 

We prove the extension by considering one cell attachment at a time. Suppose we have an extension $PW \to \C P^{n-k}$ for a subcomplex $PW$ of $PW_{n,k}$. We consider the following diagram which attaches a single cell to $PW$
$$\xymatrix{S^q\ar[d]_-{\phi}\\ 
PW\ar[r]\ar[d] &\C P^{n-k},\\
C(\phi)\ar@{-->}[ru] }  $$
and prove that the dashed arrow exists. 
This will be ensured if the composite $S^q \to \C P^{n-k}$ is trivial. 
Up to homotopy, this lifts to  $S^q\overset{\alpha}\to S^{2(n-k)+1}$, 
and it's enough to show that this lifted map is null-homotopic. We first check that $\alpha$ belongs to the stable range using \cite[Corollary 3.2]{Yam89}, which happens if $q \leq 2(n-k)p+2p-3$. Observe that 
\[
q \leq \dim(PW_{n,k}) - 1 
              = 2nk - k^2 -2  ,
\]
and that the quadratic equation in $k$ 
\[
2nk-k^2 - 2 = 2(n-k)p+2p-3 
\]
has the positive zero (over $\R$) as the second bound in the Proposition. Therefore, we are in the stable range, and further $\dim(PW_{n,k})< 2p^2-2p$ (which holds for $p> n$) guarantees that it is in fact in the image of $J$-homomorphism as in the proof of Theorem \ref{resstsplit}. We prove that the $e$-invariant of this map vanishes. 

Note that $\alpha$ induces a map $\kappa:C(\phi)\to \C P^{n-k+1}$ on mapping cones. We consider the diagram
$$\xymatrix{PW\ar[r]\ar[d] &\C P^{n-k}\ar[d]\\
                        C(\phi)\ar[r]^-\kappa\ar[d] &\C P^{n-k+1}\ar[r]\ar[d] & C\ar[d]\\
                        S^{q+1}\ar[r]^-{\Sigma \alpha} &S^{2n-2k+2}\ar[r] &C^\prime.}  $$
In this diagram the squares placed vertically come from the map between cofibre sequences and so do the squares placed horizontally. 
It suffices to show $e(\Sigma \alpha)$ vanishes. If not, there exists $\tau \in K^*(C^\prime)$ such that $ch(\tau) \not\in H^*(C^\prime;\Z_{(p)})[u^\pm]$ from \cite[Proposition 7.8]{Ada66}. Note that this implies also that the  Chern character of the image of $\tau$ in $K^*(C)$ goes outside of the image of $H^*(C;\Z_{(p)})[u^\pm] \subset H^*(C;\Q)[u^\pm]$. Our task boils down to checking that the Chern character map for $C$ takes value in the image of $H^*(C;\Z_{(p)})[u^\pm] \subset H^*(C;\Q)[u^\pm]$.

 We note that the diagram 
$$\xymatrix{PW\ar[r]\ar[d] &\C P^{n-k}\ar[d]\\
                        C(\phi)\ar[r]^-\kappa\ar[d] &\C P^{n-k+1}\ar[d] \\
                        PW_{n,k} \ar[r] & \C P^{\infty}  }  $$
commutes, where the bottom row is the classifying map of the $S^1$-bundle $S^1 \to W_{n,k}\to PW_{n,k}$.
Hence we get the following map between cofibrations,
$$\xymatrix{ C(\phi)\ar[r]^-\kappa\ar[d] &\C P^{n-k+1}\ar[d]\ar[r] & C\ar[d]^{\omega} \\
                        PW_{n,k} \ar[r] & \C P^{\infty}\ar[r] & D  } $$
with $D$ being the homotopy cofiber of $PW_{n,k}\to \C P^{\infty}$. Applying $K_{(p)}$, we obtain the diagram
$$\xymatrix{K_{(p)}^{-1}(C(\phi))\ar@{^{(}->}[r] & K_{(p)}(C)\ar[r] & K_{(p)}(\C P^{n-k+1})\ar[r]^-{\kappa^*} & K_{(p)}(C(\phi))\\
K_{(p)}^{-1}(PW_{n,k})\ar@{->>}[u]\ar@{^{(}->}[r] & K_{(p)}(D) \ar[u]\ar[r] & K_{(p)}(\C P^{\infty})\ar@{->>}[u]\ar[r] & K_{(p)}(PW_{n,k}).\ar@{->>}[u] } $$
Now  surjectivity of the two terminal vertical arrows follow from the fact that 
$C(\phi)$ is sub-complex of $PW_{n,k}$ under the minimal $p$-local $CW$ structure (Proposition \ref{p-min-cell}). 
In the proof of Theorem \ref{projstsplit}, we have checked that the image of the Chern character for $PW_{n,k}$ lies inside cohomology with $\Z_{(p)}$-coefficients, and so, the same is true for the image of $K_{(p)}^{-1}(C(\phi))$ inside $K_{(p)}(C)$. 

To complete the proof, we show that the Chern character carries $ker(\kappa^*)$ to $H^{ev}(C;\Z_{(p)})[u^\pm]$. The composition $\C P^{n-k} \longrightarrow C(\phi) \overset{\kappa}\longrightarrow \C P^{n-k+1}$ is homotopic to the inclusion,  and so, $ker(\kappa^*)$ must be equal to the $\Z_{(p)}$-module generated by $x^{n-k+1}$, where $x$ is the complex orientation of $K$-theory. Now consider the diagram obtained from \eqref{PWK}, 
$$\xymatrix{ K_{(p)}(BU(n),BU(n-k)) \ar[r]^-{f^\ast}\ar[d] & K_{(p)}(D) \ar[r]^{\omega^\ast}\ar[d] & K_{(p)}(C)\ar[d]\\
K_{(p)}(BU(n))\ar[r] & K_{(p)}(\C P^\infty)\ar[r] & K_{(p)}(\C P^{n-k+1}). }  $$
We  see that $\omega^\ast f^\ast c^K_{n-k+1}\in K_{(p)}(C)$ is mapped to $\binom{n}{n-k+1}x^{n-k+1}\in K_{(p)}(\C P^{n-k+1})$. As $\binom{n}{n-k+1}$ is a unit in $\Z_{(p)}$, it suffices to show that $ch(\omega^\ast f^\ast c^K_{n-k+1})$ lies in $H^{ev}(\C P^{n-k+1};\Z_{(p)})[v^\pm]$.  We also notice that using computations analogous to Lemma \ref{chgammaexp}
\begin{align*}
ch(\omega^\ast f^\ast c^K_{n-k+1}) &= \omega^\ast f^\ast (ch(c^K_{n-k+1}))\\
&= \omega^\ast f^\ast (c_{n-k+1}\cdot P_{n-k+1}+ \cdots + c_n \cdot P_n) \text{ [where } P_j \in \Q[[c_1, \cdots, c_n]] \text{ ]}\\
&= \omega^\ast f^\ast c_{n-k+1} \cdot Q_{n-k+1} +\cdots + \omega^\ast f^\ast c_n \cdot Q_n \text{ [where } Q_j \in \Q[x]/(x^{n-k+2}) \text{ ]}
\end{align*}
This computation shows that the maximum degree of the homogeneous parts of $ch(\omega^\ast f^\ast c^K_{n-k+1})$ is 
\[
|c_n| + |x^{n-k+1}|= 2n + 2(n-k+1). \]
As $\kappa: C(\phi) \to \C P^{n-k+1}$ is a $(2n-2k+1)$-equivalence,  $ C$ is  $(2n-2k)$-connected. 
Therefore, as $p>n+1$, we apply \cite[Theorem 1]{Ada61} to deduce the result.                                                                          
\end{proof}

Proposition \ref{retcp} constructs for us the map from $PW_{n,k} \to \C P^{n-k}$. The remaining factors in the product decomposition of $Y_{n,k}$ are spherical and maps are constructed via a connectivity argument. 
\begin{prop}\label{sphfactor}
If $p>n$ and $k \leq \text{min}\{n,(p+n)- \sqrt{p^2+n^2-4p+2}\}$, then for $1 \leq r \leq k-1$, there is a map from $PW_{n,k}$ to $S^{2(n-k)+2r+1}$ which pulls back the standard $\Z_{(p)}$-cohomology generator of $ S^{2(n-k+r)+1}$ to the class $\gamma_{n-k+r}\in H^*(PW_{n,k};\Z_{(p)})$. 
\end{prop}
\begin{proof}
From the equivalence of $\Sigma^{\infty} PW_{n,k}$ and $\Sigma^{\infty} Y_{n,k}$ we get maps  
\[\nu_r: \Sigma^{\infty} PW_{n,k} \to \Sigma^{\infty}S^{2(n-k+r)+1}\]
 for $1 \leq r \leq k-1$, which satisfy $\nu_r^\ast (\epsilon_{2(n-k+r)+1}) = \gamma_{n-k+r}$. By the usual $\Sigma^\infty-\Omega^\infty$ adjunction, we obtain a map of spaces 
\[\tilde{\nu}_r : PW_{n,k} \to QS^{2(n-k+r)+1}.\]
 Now as we are in the $p$-local category, the fibre $F(2(n-k+r)+1)$ of the natural map $S^{2(n-k+r)+1}\to  QS^{2(n-k+r)+1}$ is $(2p(n-k+r)+2p-4)$-connected \cite[corollary 3.2]{Yam89}. The given bounds on $k$ imply that $\dim(PW_{n,k})\leq 2p(n-k+r)+2p-4$. Hence, the map $\tilde{\nu_r}$ lifts to $S^{2(n-k+r)+1}$, and we are done. 
%
\end{proof}

We let $M(n,p)= p+n-\sqrt{p^2+n^2 - 4p +2}$, and summarize the results of Propositions \ref{retcp} and \ref{sphfactor} in the following theorem. 
\begin{thm} \label{unsplit}
Let $p>n+1$ and $k\leq \text{min}(n,M(n,p))$. Then, in the $p$-local category 
\[
PW_{n,k} \simeq \C P^{n-k} \times S^{2n-2k+3} \times \cdots \times S^{2n-1}.\] 
\end{thm} 

Observe that the bound on $k$ holds whenever $k<n/2$. As in \S \ref{decvlarge}, the product decomposition of the projective Stiefel manifold implies an $S^1$-equivariant decomposition of the complex Stiefel manifold. 

\begin{thm}\label{eqstief}
Let $p>n+1$, and $k\leq \text{min}(n,M(n,p))$. Then, we have an equivalence of $S^1$-spaces  
\[
{W_{n,k}}_{(p)} \simeq \Big[ S^{2(n-k)+1} \times S^{2n-2k+3} \times \cdots \times S^{2n-1}\Big]_{(p)}.\] 
\end{thm} 
\begin{proof}
 As both sides of the equivalence possess a free $S^1$-action, it suffices to exhibit an $S^1$-equivariant map which is a weak equivalence. We use the map 
\[ \phi_k : PW_{n,k} \to \C P^{n-k}\]
from Proposition \ref{retcp} and pullback the circle bundle $q: S^{2(n-k)+1} \to \C P^{n-k}$ via $\phi_k$. As $\phi_k$ is a lift of the classifying map of the circle bundle $W_{n,k}\to PW_{n,k}$, we have a commutative square 
\[
\xymatrix{ W_{n,k} \ar[r]^{\hat{\phi}_k} \ar[d] & S^{2n-2k+1} \ar[d]\\ 
 PW_{n,k} \ar[r]  & \C P^{n-k}. } \] 
We now observe that $\hat{\phi}_k$ is $S^1$-equivariant as it fits in a pullback diagram of $S^1$-bundles. We now form the equivariant map 
\[ 
W_{n,k} \to S^{2n-2k+1} \times W_{n,k-1}.
\]
Now $k-1 \leq M(n,p)$, and we proceed by induction assuming that $W_{n,k-1}$ supports the splitting stated in the theorem, to deduce the result.
\end{proof}

There are also analogous splittings for the spaces $W_{n,k;m}$ and $P_\ell W_{n,k}$. 

\begin{thm}\label{splitquot}
Let $p>n+1$, and $k\leq \text{min}(n,M(n,p))$. Then, we have equivalences
\begin{align*}
\Big(W_{n,k;m}\Big)_{(p)} &\simeq \Big[ L_m(2n-2k+1) \times S^{2n-2k+3} \times \cdots \times S^{2n-1}\Big]_{(p)}~~ \mbox{ if } p\mid m,\\
\Big(P_\ell W_{n,k}\Big)_{(p)} & \simeq \Big[ \C P^{n-k} \times S^{2n-2k+3} \times \cdots \times S^{2n-1} \Big]_{(p)}.
\end{align*}
\end{thm}

\begin{proof}
The $S^1$-equivariant map $W_{n,k}\to S^{2n-2k+1}$ in Theorem \ref{eqstief} yields the map 
\[ W_{n,k;m} \to L_m(2n-2k+1) \]
on $C_m$-orbits. The other factors receive maps via the composite 
\[ W_{n,k;m} \to PW_{n,k} \to S^{2n-2k+2r+1} \]
where the latter map is defined by the equivalence of Theorem \ref{unsplit}. The product of these two maps imply the first $p$-local equivalence in the theorem. For the second equivalence, one applies the stable decomposition for $P_\ell W_{n,k}$ outlined in Remark \ref{splplwnk}. After that, observe that Proposition \ref{sphfactor} works verbatim for $P_\ell W_{n,k}$ and Proposition \ref{retcp} works analogously by using the pair $(BU(k), Gr_k(\C^n))$ instead of $(BU(n), BU(n-k))$. 
\end{proof}

\end{document}